\date{}

\date{}

\title{On invariant generating sets for the cycle space}
\author{\'Ad\'am  Tim\'ar }

\documentclass[11pt,naturalnames]{article}
\renewcommand\footnotemark{}
\usepackage[ruled,vlined]{algorithm2e}
\usepackage{amsmath}
\usepackage{amsthm}
\usepackage{amsfonts}
\usepackage{graphicx}
\newif\ifhyper\IfFileExists{hyperref.sty}{\hypertrue}{\hyperfalse}
\ifhyper\usepackage{hyperref}\fi
\usepackage{enumitem}
\usepackage{subfig}
\usepackage{caption}
\usepackage{keyval}
\usepackage[margin=1in]{geometry}
\usepackage{setspace}

\usepackage[displaymath, mathlines,pagewise]{lineno}

\theoremstyle{definition}

\newtheorem{theorem}{Theorem}

\newtheorem{lemma}[theorem]{Lemma}

\newtheorem{remark}[theorem]{Remark}

\def \proof {{ \medbreak \noindent {\bf Proof.} }}
\def\proofof#1{{ \medbreak \noindent {\bf Proof of #1.} }}
\def\proofcont#1{{ \medbreak \noindent {\bf Proof of #1, continued.} }}
\def\supp{{\rm supp}}
\def\max{{\rm max}}
\def\min{{\rm min}}
\def\dist{{\rm dist}}
\def\Aut{{\rm Aut}}
\def\id{{\rm id}}
\def\Stab{{\rm Stab}}
\def\span{{\rm span}}

\begin{document}
\maketitle
\let\thefootnote\relax\footnotetext{\footnotesize{Partially supported by Icelandic Research Fund grant number 239736-051 and the ERC Consolidator Grant 772466 ``NOISE''.}}

\bigskip

\def\eref#1{(\ref{#1})}
\newcommand{\Prob} {{\bf P}}
\newcommand{\calC}{\mathcal{C}}
\newcommand{\calP}{\mathcal{P}}
\newcommand{\calQ}{\mathcal{Q}}
\newcommand{\calO}{\mathcal{O}}
\newcommand{\LL}{\mathcal{L}}
\newcommand{\Z}{\mathbb{Z}}
\newcommand{\N}{\mathbb{N}}
\newcommand{\HH}{\mathbb{H}}
\newcommand{\Rr}{\mathbb{R}^3}
\newcommand{\h}{\mathcal{H}}
\def\diam{\mathrm{diam}}
\def\length{\mathrm{length}}
\def\ev#1{\mathcal{#1}}
\def\Isom{{\rm Isom}}
\def\Re{{\rm Re}}
\def \eps {\epsilon}
\def \P {{\Bbb P}}
\def \E {{\Bbb E}}
\def \proof {{ \medbreak \noindent {\bf Proof.} }}
\def\proofof#1{{ \medbreak \noindent {\bf Proof of #1.} }}
\def\proofcont#1{{ \medbreak \noindent {\bf Proof of #1, continued.} }}
\def\supp{{\rm supp}}
\def\max{{\rm max}}
\def\min{{\rm min}}
\def\dist{{\rm dist}}
\def\Aut{{\rm Aut}}
\def\id{{\rm id}}
\def\Stab{{\rm Stab}}
\def\T{{\cal T}}
\def\B{{\cal B}}

\newcommand{\lra}{\leftrightarrow}
\newcommand{\xlra}{\xleftrightarrow}
\newcommand{\xnlra}{\xnleftrightarrow}
\newcommand{\pc}{{p_c}}
\newcommand{\pt}{{p_T}}
\newcommand{\ptk}{{\hat{p}_T}}
\newcommand{\pl}{{\tilde{p}_c}}
\newcommand{\pe}{{\hat{p}_c}}
\newcommand{\pr}{\mathrm{\mathbb{P}}}
\newcommand{\pp}{\mu}
\newcommand{\ex}{\mathrm{\mathbb{E}}}
\newcommand{\ee}{\mathrm{\overline{\mathbb{E}}}}

\newcommand{\om}{{\omega}}
\newcommand{\ebd}{\partial_E}
\newcommand{\ivbd}{\partial_V^\mathrm{in}}
\newcommand{\ovbd}{\partial_V^\mathrm{out}}
\newcommand{\q}{q}
\newcommand{\RR}{\mathcal{R}}

\newcommand{\CC}{\Pi}
\newcommand{\BB}{\Pi}

\newcommand{\A}{\mathcal{A}}
\newcommand{\cc}{\mathbf{c}}
\newcommand{\pa}{{PA}}
\newcommand{\degi}{\deg^{in}}
\newcommand{\dego}{\deg^{out}}
\def\Pn{{\bf P}_n}

\newcommand{\R}{\mathbb R}
\newcommand{\F}{F}
\newcommand{\FF}{\mathfrak{F}}
\newcommand{\Can}{\rm Can}
\newcommand{\Vol}{\mathrm Vol}

\def\UST{{\rm UST}}
\def\Gstar{{{\cal G}_{*}}}
\def\Gstarstar{{{\cal G}_{**}}}
\def\Gstarplus{{{\cal G}_{*}^\frown}}
\def\Rel{{\cal R}}
\def\Comp{{\rm Comp}}
\def\calG{{\cal G}}
\def\calE{{\cal E}}
\def\omps{{\omega_\delta^\eps}}
\def\FKI{{\phi^f_{G,p,0}}}
\def\FKIn{{\phi_{G_n,p,0}}}
\def\FKIf{{\phi^f_{G,p,0}}}
\def\FKIw{{\phi^w_{G,p,0}}}
\def\loop{{\Prob_{G,x}}}
\def\loopn{{\Prob_{G_n,x}}}

\begin{abstract}
Consider a unimodular random graph, or just a finitely generated Cayley graph. When does its cycle space have an invariant random generating set of cycles such that every edge is contained in finitely many of the cycles?
Generating the free Loop $O(1)$ model as a factor of iid is closely connected to having such a generating set for FK-Ising cluster. We show that geodesic cycles do not always form such a generating set, by showing for a parameter regime of the FK-Ising model on the lamplighter group 
the origin is contained in infinitely many geodesic cycles. This answers a question by Angel, Ray and Spinka. Then we take a look at how the property of having an invariant locally finite generating set for the cycle space is preserved by Bernoulli percolation, and apply it to the problem of factor of iid presentations of the free Loop $O(1)$ model.
\end{abstract}

\bigskip

\section{Introduction}

For a finite graph $G$ the space $\calE$ of even subgraphs is the space of subgraphs where every degree is even, with modulo 2 addition. The set of (simple) cycles generates it. 
The {\it random/uniform even subgraph} is a uniformly chosen random element of $\calE$.
If $G$ is infinite, and a sequence of finite graphs locally converges to it, one can take the wired or free boundary conditions, and take the weak limit of the random even subgraph. The resulting random element of $\calE$ is called the {\it wired/free uniform even subgraph} of $G$, as defined in \cite{GJ}. One can prove that this is invariant under the automorphisms of $G$.

When is it possible to find a set $\calC$ of cycles that generates $\calE$ and such that $\calC$ is {\it locally finite}? This latter means that every edge is contained in finitely many elements of $\calC$. Grimmett and Janson \cite{GJ} proved that there is always such a generating set. But suppose that we want to find a random generating set, and such that its distribution is invariant under the automorphisms of $G$. This is a natural question already when $G$ is a Cayley graph (or a unimodular transitive graph). For the wired case (where one has to include biinfinite paths among the cycles) there is a simple negative answer, as can be shown for the regular tree. Therefore we are only interested in the free case; also because the question that originally motivated us is only open in this setup. We present this question next.

A random subgraph of a transitive graph $G$ is called a {\it factor of iid}, if one can encode it as an equivariant and measurable function from the space $[0,1]^{V(G)}$ of iid labellings of the vertices of $G$ by Lebesque[0,1] labels.
The Loop $O(1)$ model is in the focus of recent work by Angel, Ray and Spinka \cite{ARS}, inspired by \cite{H}. The authors pursue a factor of iid construction through the known factor of iid construction of the FK-Ising model (\cite{HJL,HSr}). 
While they gave a full solution for the wired model, the free case was only partially solved, so we restrict ourselves to this latter now.
The problem of constructing the free Loop $O(1)$ model as a factor of iid is equivalent to finding a random even subgraph as a factor of iid.
One way to achieve this would be to find a locally finite invariant generating sets of cycles, and moreover, find it as a factor of iid. Angel, Ray and Spinka examined the family of {\it geodesic cycles}: a cycle is {\it geodesic} in graph $G$ if for every two points of the cycle at least one of the arcs between them is geodesic in $G$. The set of all geodesic cycles is a factor of iid and is always a generating set for the cycle space. Thus, if the set of geodesic cycles in the FK-Ising configuration happens to be locally finite, then the above problem is settled. 
This has led Angel, Ray and Spinka ask if this is always the case \cite{ARS}. We give a negative answer to their question.

\begin{theorem}\label{DL}
Let $G={{\rm DL}}(2,2)$ be the Diestel-Leader representation of the lamplighter group with some fixed vertex $o$. There is a $\tilde p<1$ such that if $p>\tilde p$, then with positive probability $\omega_p$ contains infinitely many geodesic cycles through $o$. Here $\omega_p$ may be the Bernoulli($p$) bond percolation configuration or the FK-Ising configuration $\FKI$.
\end{theorem}

When $\omega_p$ is the Bernoulli($p$) bond percolation, one can choose $\tilde p=\frac{1}{\sqrt{2}}$ in the theorem. Then for $\omega_p=\FKI$, the choice $\tilde p=\frac{2}{\sqrt{2}+1}$ works, by Lemma \ref{stochdom}.

\begin{remark}\label{nemfifi}
If we take $G$ to be the direct product of the 3-regular tree and ${{\rm DL}}(2,2)$, we get a {\it nonamenable} transitive graph where the conclusion of Theorem \ref{DL} holds, using Theorem \ref{DL} and Lemma \ref{stochdom}. This is an example where it is clear that none of the conditions in \cite{ARS} imply that the
free Loop $O(1)$ model would be a factor of iid.
\end{remark}

One way to find a factor of iid locally finite set of cycles generating the cycle space is through geodesic cycles, another one is obtained when a one-ended factor of iid spanning tree exists, in which case one can take all cycles that have a single edge outside of the tree. As Remark \ref{nemfifi} shows, there are examples where neither of these approaches can be successful. (Note that a nonamenable unimodular transitive graph cannot have an invariant one-ended spanning tree.) So one may wonder about other approaches, and one reasonable question is whether percolation can change the existence of a locally finite generating sets of cycles. These are the topics of our last section. 
In particular, let $G$ be the Cayley graph of a finitely presented group, in which the cycles given by all conjugates of a finite set of relators defines a locally finite generating set for $\calE$. Let $\omega$ be a supercritical percolation configuration. Is there a locally finite unimodular/factor of iid set of generators for the cycle space of $\omega$?

\begin{theorem}\label{percol}
If $G$ is the Cayley graph of a finitely presented group, then there exists a $\bar p<1$ such that $G(p)$ has an invariant feasible generating set whenever $\bar p<p<1$, and this can be chosen as a factor of iid. The same holds for $\FKI$. 
\end{theorem}

The existence of such a factor of iid generating set would imply right away that the Loop $O(1)$ model on Cayley graphs of hyperbolic groups is a factor of iid. (Also see Question 6.3 and related remarks in \cite{ARS}.) In Theorem \ref{percol} we answer the above question for percolation probability close enough to 1, which implies the following theorem.

\begin{theorem}\label{finitely_pres}
Let $G$ be the Cayley graph corresponding to a finite presentation of some finitely presented group, such as a hyperbolic group.
Then there is some $\bar x<1$ such that the free Loop $O(1)$ model of parameter $x$ and 0 on $G$ is a factor of iid, for every $\bar x<x\leq 1$.
\end{theorem}
The proof of this theorem readily extends to unimodular random graphs whose cycle space is generated by cycles of length at most $k$ with some $k$. 

After providing the necessary definitions in Section \ref{definitions}, including those that have already come up, we prove Theorem \ref{DL} in Section \ref{geodesic}. Some further thoughts on invariant/factor of iid generating sets of the cycle space are in Section \ref{concluding}, together with a proof of Theorem \ref{percol} and its corollary, Theorem \ref{finitely_pres}.

\section{Definitions, preliminaries}\label{definitions}

Let ${\cal G}_{*}$ be the set of locally finite graphs $(G,o)$ with a root vertex up to isomorphisms fixing the root. Let ${\cal G}_{**}$ be the set of locally finite graphs $(G,o,x)$ with an ordered pair of vertices, up to isomorphisms preserving the ordered pair. 
A probability measure $\mu$ on ${\cal G}_{*}$ is called a {\it unimodular random graph (URG)} if the following {\it Mass Transport Principle (MTP)} holds:
\begin{equation}
\int \sum_{x\in V(G)} f(G,o,x) d\mu ((G,o))=\int \sum_{x\in V(G)} f(G,x,o) d\mu ((G,o))
\label{eq:MTP}
\end{equation}
for every Borel function $f:{\cal G}_{**}\to [0,\infty]$. 
Sometimes we will simply refer to the random rooted graph $(G,o)$ as a URG. 
An important source of URG's is to take a random subgraph of a Cayley graphs (or more generally, a unimodular transitive graph) whose distribution is invariant under the isomorphisms of the graph (or: jointly unimodular with $G$), and look at the component of a fixed origin. We are interested in generating the Loop $O(1)$ model as a uniform even subgraph of the $\FKI$ on $G$. Hence URG's seem to be the natural choice of generality when we study the cycle space of the $\FKI$ clusters of a Cayley graph.
We refer the reader to \cite{AL} or \cite{Cu} for more background on URG's.

Let $X$ be some separable metric space. 
A Borel measurable function $f : {\cal G}_* \to X$ is called a {\it factor map}.
Let G be a random graph, and
$\lambda : V (G)\to [0, 1]$ be iid Lebesgue[0,1] random labels
on its vertices, and denote by $G(\lambda)$ the resulting random labeled graph. If $f$ is a factor map, then the set of random variables 
$\bigl(f((G(\lambda),v))\bigr)_{v\in V(G)}$
is
a {\it graph factor of iid} or just {\it factor of iid}. 
That is, a factor of iid is a new labelling of $V(G)\cup E(G)$ with elements from some metric space $X$ with the property that
for every $\eps>0$ there is an $r$ such that with probability at least $1-\eps$, from the Lebesgue labels in the $r$-neighborhood of $x\in V(G)\cup E(G)$ one can tell its new label, up to an error of $<\eps$ in $X$.

A unimodular probability measure $\mu$ on ${\cal G}_{*}$ is called {\it amenable} if 
there exist Borel functions $\lambda_n : {\cal G}_{**}\to [0,1]$ such that  
$$||\lambda_n(G,o,.)||_1=1
$$
and
$$\lim_{n\to\infty}\sum_{x\sim o } ||\lambda_n (G,o,.)-\lambda_n (G,x,.)||_1=0
$$
for $\mu$ almost every $(G,o)$.


Now we define the {\it Fortuin-Kasteleyn Ising model} (or {\it FK-Ising model}). For us it will be enough to define the special case when the external field is 0 (see \cite{PGG} for the more general case and for the larger class of random cluster models). Given a finite graph $G=(V,E)$ and parameter value $p\in [0,1]$, define
\begin{equation}\label{fki}
{\phi_{G,p,0}} (\omega)=\frac{1}{Z_{G,p,0}} \Bigl(\frac{p}{1-p}  \Bigr)^{|\{e\in E:\omega(e)=1 \}|}2^{k(\omega)},
\end{equation}
where $k(\omega)$ is the number of connected components of $\omega$ and $Z_{G,p,0}$ is a normalizing constant.
If $G$ is an infinite graph and $G_n$ is an increasing sequence of finite graphs whose union is $G$, then $\FKIn$ has a weak limit $\FKI$, and it is the same for any such sequence $G_n$. We call this limit the {\it free FK-Ising measure}. One can analogously define a wired FK-Ising measure, but this will not be needed for the present paper.
It has been proved that the FK-Ising models are factors of iid, \cite{HJL,HSr} (see also the explanation in \cite{ARS}).


An {\it even subgraph} of $G$ is a spanning subgraph (meaning a subgraph that contains all vertices), where every vertex has an even degree. If $G$ is finite, we will call the probability measure 
\begin{equation}
\loop (\omega)=\frac{1}{Z_{G,x}} x^{|\{e\in E: e\in\omega \}|}
\end{equation}
on the set $\calE=\{\omega \text{ an even subgraph of }G\}$
the {\it Loop $O(1)$ model of parameter} $x$. Here $Z_{G,x}$ is the normalizing constant. Setting $x=1$, we get the {\it uniform even subgraph} of $G$.
As before, we mention that the above definitions can be made more general, but for the questions that we are interested in, the present definitions are sufficient (with parameters $y=0$ and $B=\emptyset$, in the notation of \cite{ARS}). Now, if $G$ is infinite, one can take an increasing sequence of finite graphs $G_n$ whose union is $G$, and $\loopn$ will weakly converge to a probability measure $\loop$ on the even subgraphs of $G$, \cite{ADS}. We call this measure the {\it free Loop $O(1)$ measure} (of parameter $x$) on $G$.

For a given graph $G=(V,E)$, consider the vector space over the 2-element field of all spanning subgraphs. So addition of two subgraphs is defined as modulo 2 addition of their edges, that is, the sum of the two sets is their symmetric difference. The set of even subgraphs of $G$ forms a subspace of this vector space. If $G$ is finite, then this subspace is generated by the set of finite cycles (which may not be true anymore when $G$ is infinite, e.g. a biinfinite path), and so one can identify the set of even subgraphs with the cycle space (to be defined soon for an infinite $G$). 
Now for $G$ infinite, consider the set of all finite sums of cycles, their closure in the topology of pointwise convergence, and call the resulting subspace 
the {\it (free) cycle space} of $G$. If $G_n$ is an increasing sequence of finite graphs whose union is $G$, take a uniform random element of the cycle space of $G_n$. This uniform even subgraph weakly converges to a random element of the free cycle space of $G$, called the {\it free uniform even subgraph} of $G$.

The connection between the free Loop $O(1)$ model as factor of iid and locally finite invariant generating sets of cycle spaces is summarized in the next theorem, as a concequence of previous results.
Namely, the Loop $O(1)$ model can be attained as a uniform even subgraph of the (free) FK-Ising configuration by \cite{ADTW}, see Proposition 2.6 and Remark 2.7 in \cite{ARS}. From \cite{HJL} and \cite{HSr}
we know that $\phi_{G,p,0}$ is a factor of iid, moreover, it is finitary. Finally, any factor of iid locally finite generating set of the cycle space of a graph can be used to construct a factor of iid representation of the uniform even subgraph of that graph (Proposition 3.12 in \cite{ARS}. These combine to the following.

\begin{theorem}[\cite{ADTW},\cite{HJL},\cite{HSr}] 
(1) The free random even subgraph of the free FK-Ising $\phi_{G,p,0}$ sample
has the same distribution as the free Loop $O(1)$ model of parameter $x=\frac{p}{2-p}$.\\
(2) The free Loop $O(1)$ is a factor of iid (respectively: finitary factor of iid) if there is a factor of iid (respectively: finitary factor of iid) locally finite generating set for the cycle space of $\phi_{G,p,0}$.
\end{theorem}


\section{Infinitely many geodesic cycles though a point}\label{geodesic}
In \cite{ARS}, Angel, Ray and Spinka ask the following question (Question 6.2). Let $G$ be a unimodular random graph and consider the FK-Ising model $\FKI$ on $G$ with some supercritical parameter $p<1$. Is it possible that $\FKI$ contains infinitely many geodesic cycles through the root? We provide a positive answer, taking a Cayley graph of the lamplighter group as an example, with $p$ large enough.
The motivation for the question was coming from Theorem 1.2 in \cite{ARS}. There the authors prove that the free Loop $O(1)$ model is a factor of iid if $y>0$, or if $G$ is invariantly amenable, or if the following (``third'') condition holds. For the sample $\omega$ of the free FK-Ising measure on $G$ with parameter $p=\frac{2x}{1+x}$, there are only finitely many geodesic cycles through any vertex almost surely.
The authors say,
``Regarding the third condition, we do not know whether it is satisfied in general
for $x<1$.'' In light of our Theorem \ref{DL}, the condition is not satisfied in general.
A negative answer to Question 6.2 in \cite{ARS} would have implied that their factor of iid construction works for every $G$. Our example proves the opposite, and it also provides an example where the results of \cite{ARS} are known to not imply that the free FK-Ising model is a factor of iid; see Remark \ref{nemfifi}.

In Section \ref{concluding} we will propose replacements to the local finiteness of the set of geodesic cycles. 

We will need the following well known consequence of the FKG-inequality for the FK model. See e.g. Subsection 13.1 of \cite{PGG}.

\begin{lemma}\label{stochdom}
The FK-Ising measure $\FKI$ stochastically dominates Bernoulli($\tilde p$) bond percolation whenever $\tilde p\leq \frac{p}{2-p}$. 
\end{lemma}

Let us recall the {\it Diestel-Leader (horocyclic product) representation of the lamplighter group}. Consider two 3-regular infinite trees, $T_1,T_2$, together with graph homomorphisms $h_1$ of $T_1$ and $h_2$ of $T_2$ into $\Z$. The map $h_1$ has the property that for every $x$ in $V(T_1)$ the value of $h_1$ is $h_1(x)-1$ in exactly two neighbors of $x$ (call their set $N$ for the moment), and it is $h_1(x)+1$ in exactly one neighbor. If a vertex $y$ is separated from $x$ by $N$, then we call $y$ a {\it descendant} of $x$, and call $x$ an {\it anchestor} of $y$.
We can think of $T_1$ as a family tree where everybody has one parent and two children, and $h_1$ tells the generation that a vertex belongs to (generations numbered backwards in time). The mapping $h_2$ has the property that for every $x$ in $V(T_2)$ the value of $h_2$ is $h_2(x)+1$ in exactly two neighbors of $x$, and it is $h_2(x)-1$ in exactly one neighbor. 
We will use the terms descendant and ancestor for $T_2$, so every vertex has one neighbor that is its ancestor and two neighbors that are its descendants.
Now let $G={{\rm DL}}(2,2)$ be defined on vertex set $\{(v_1,v_2)\in V(T_1)\times V(T_2), h_1(v_1)+h_2(v_2)=0\}$, with $(x_1,x_2)$ and $(y_1,y_2)$ adjacent if and only if $\{x_1,y_1\}\in E(T_1)$ and $\{x_2,y_2\}\in E(T_2)$. 
Given a vertex $(v_1,v_2)\in V(G)$, we will refer to $v_1\in V(T_1)$ as its first coordinate, and $v_2\in V(T_2)$ as its first coordinate, Similarly, any edge $e\in E(G)$ can be written in the form $(e_1,e_2)$, where, $e_1\in E(T_1)$ and $e_2\in E(T_2)$. Then we will call $e_1$ (respectively $e_2$) the first (respectively second) {\it coordinate of} $e$.
By a {\it path} we always mean a path that is not self-intersecting, that is, no repetition of vertices is allowed.

\begin{proofof}{Theorem \ref{DL}}
By Lemma \ref{stochdom} it is enough to prove the claim for Bernoulli percolation. Fix $o=(o_1,o_2)\in V(G)$. We show that for every $p>1/\sqrt{2}$, the Bernoulli($p$) bond percolation $\omega_p$ contains infinitely many geodesic cycles through $o$ with positive probability. 

For a binary tree with root $\rho$ (that is, a tree where every degree is 3, except for $\rho$, which has degree 2), denote by $\theta(p)$ the probability that $\rho$ is in an infinite cluster by Bernoulli($p$) bond percolation. 

Next we give names to a few vertices and subsets of vertices, which is all illustrated on Figure \ref{names}. To facilitate reading, vertices in $T_1$ will always have index 1 (even when there is no corresponding point of index 2...), and vertices in $T_2$ will always have index 2. When a vertex has no index then it is a vertex of $G$.
Consider the subtree induced by all descendants of $o_1$ in $T_1$ up to distance $n$ from $o_1$, and call this subtree $T_{1,n}$.
Descendant vertices at distance $n$ from $o_1$ belong to one of two classes, depending on which child of $o_1$ is their ancestor. Call these two classes $L_1$ and $L_1'$.
Denote by $\hat o_2$ the unique ancestor of $o_2$ at distance $n$ from it in $T_2$. Call the subtree of $T_2$ induced by the set of descendants of $\hat o_2$ in $T_2$ up to distance $n$ from $\hat o_2$ as $T_{2,n}$. Let $L_2$ be the set of descendant vertices of $\hat o_2$ in $T_2$ at distance $n$ from $\hat o_2$ that are separated from $o_2$ by $\hat o_2$ (so they belong to the other descendant subtree from $\hat o_2$). 

\begin{figure}[h]
\vspace{0.1in}
\begin{center}
\includegraphics[keepaspectratio,scale=1.]{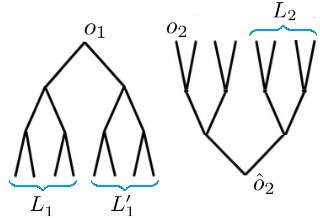}
\caption{The two subtrees $T_{1,n}$ and $T_{2,n}$.
We called $T_o$ the edges of $\omega_p$ whose first coordinate is in the left tree and second coordinate is in the $o_2$-$\hat o_2$ path on the right. It is distributed as Bernoulli($p$) percolation on the left tree, and thus will have some vertices ${\ell_1}\in L_1$ and ${\ell'_1}\in L_1'$ in the cluster of $o_1$ with probability at least $p^2\theta(p)^2$. Then we look at the $T_{{\ell_1}}\subset\omega_p$ of edges whose first coordinate is in the $o_1$-${\ell_1}$ path on the left and second coordinate is in the right tree. Define $T_{{\ell'_1}}$ similarly. With uniformly positive probability, the two independent Bernoulli($p$) configurations that arise if we project $T_{{\ell_1}}$ and $T_{{\ell'_1}}$ to the right tree, both contain some element $(o_1,x_2)$ with $x_2\in L_2$. A geodesic cycle in $\omega_p$ is then given by the open paths between the consecutive pairs $o$, $({\ell_1},\hat o_2)$, $(o_1,x_2)$, $({\ell'_1},\hat o_2)$, $o$.
}
\label{names}
\end{center}
\end{figure}

Consider the subtree $T_o$ of $G$ consisting of edges whose first coordinate is in
$T_{1,n}$ and second coordinate is in the $o_2$-$\hat o_2$ path of $T_2$. Let $\omega(L_1)$ consist of those such $v_1\in L_1$ that are in the same component as $o_1$ when we restrict $\omega_p$ to this subtree. Define $\omega(L_1')$ similarly, replacing $L_1$ by $L_1'$ in the above definition. 
The sets $\omega(L_1)$ and $\omega(L_1')$ are the $n$'th generations of supercritical random trees given by Bernoulli($p$) bond percolation on two (edge-disjoint) binary trees (with an extra edge to $o_1$), hence the probability that both $\omega(L_1)$ and $\omega(L_1')$ are nonempty is at least $p^2\theta(p)^2$. Condition on the event $E$ that $\omega(L_1)$ and $\omega(L_1')$ are nonempty, and fix some ${\ell_1}\in L_1$ and ${\ell'_1}\in L_1'$. 

Now we look at the subtree $T_{{\ell_1}}$ (respectively $T_{{\ell'_1}}$) of $G$ consisting of edges whose second coordinate is in
$T_{2,n}$ and first coordinate is in the ${\ell_1}$-$o_1$ path (respectively ${\ell'_1}$-$o_1$ path) of $T_1$. Note that $T_{{\ell_1}}$ and $T_{{\ell'_1}}$ are edge-disjoint, because the ${\ell_1}$-$o_1$ path and ${\ell'_1}$-$o_1$ path are edge-disjoint in $T_1$.
Define $X$ (respectively $X'$) to be the set of $v_2\in L_2$ such that there is a path in $\omega_p$ from ${\ell_1}$ (respectively ${\ell'_1}$) to $(o_1,v_2)$ in $T_{{\ell_1}}$ (respectively $T_{{\ell'_1}}$). 
By the edge-disjointness of $T_{{\ell_1}}$ and $T_{{\ell'_1}}$, $X$ and $X'$
are independent. Furthermore, both $X$ and $X'$ are the $n$'th generation of a Bernoulli($p$) bond percolation ($\omega_p$) on a binary tree, therefore 
\begin{equation}\label{nemures}
\Prob (X\cap X'\not=\emptyset|E)\geq \Prob (X\cap X'\not=\emptyset)\geq \theta(p^2),
\end{equation}
where for the first inequality we were using the FKG inequality.

We have obtained that the probability that ${\ell_1}$ and ${\ell'_1}$ exist, and that some $x_2\in
X\cap X'$ exists is at least $p^2\theta(p)^2\theta(p^2)>0$
for every $n$.
On this event, consider the union of the following pairwise edge-disjoint four paths:
\begin{itemize}
\item the path from from $o$ to $({\ell_1},\hat o_2)$ in $T_o$;
\item the path from $({\ell_1},\hat o_2)$ to $(o_1,x_2)$ in $T_{{\ell_1}}$;
\item the path from $(o_1,x_2)$ to $({\ell'_1},\hat o_2)$ in $T_{{\ell_1'}}$;
\item the path from $({\ell'_1},\hat o_2)$ to $o$ in $T_o$.
\end{itemize}
The resulting union is a cycle $C$ of length $4n$ in $\omega_p$ through $o$. we show next that this cycle is geodesic. Say that two vertices are at the same level if $h_1$ of their first coordinates is the same (that is, their coordinates belong to the same generation in $T_1$ and $T_2$). 
Note that each of the four paths is monotonous in the sense that the path keeps going up the levels or down the levels. From this property it follows that if two points of $C$ are at distance $2n$ within $C$ (that is, if the two points are {\it diagonal}), then they have to be on the same level. Now if $x$ and $y$ are two arbitrary diagonal points in $C$, to find a shortest path between them in $G$ one would first have to go from the level of $x$ up the level of $o$ (to be able to swith between the two halves of the descendant subtree of $o_1$ where their first coordinates belong to), then down to the level of $\ell_1$ and $\ell_1'$, and then back to the level of $y$ (which is the level of $x$). This means a minimum number of $2n$ steps. So the $G$-distance between any two points of $C$ at $C$-distance $2n$ is $2n$. Now, if $C$ had two points $x$ and $y$ whose $C$-distance was strictly larger than their $G$-distance, then by replacing $y$ by the point $y'$ that is diagonal to $x$, one would get 
$$\dist_G (x,y')\leq \dist_G (x,y)+\dist_G(y,y')<\dist_C (x,y)+\dist_G(y,y')\leq \dist_C (x,y)+\dist_G(y,y')=2n,$$
contradicting that diagonal points have $G$-distance $2n$. This shows that $C$ is geodesic.

We showed that for every $n$, with probability bounded below by some uniform positive constant, there is a geodesic cycle of length $4n$ through $o$ in $\omega_p$. The Borel-Cantelli lemma then implies the theorem.
\qed
\end{proofof}




\section{Locally finite generating sets and percolation}\label{concluding}

Geodesic cycles became of interest in \cite{ARS} because by flipping a coin on each of them, and making an edge open if it is contained in an odd number of the cycles, we end up with the uniform even subgraph. For this to make sense, one of course needs to make sure that no edge is contained in an infinite number of these cycles, hence the importance of Question 6.2 in \cite{ARS}. 
But the only characteristic of geodesic cycles that is needed for this application in \cite{ARS} is that the set of all geodesic cycles is automatically a factor of iid generating set for the cycle space. (This is so because a cycle that is not geodesic can be written as the sum of two shorter cycles.) Any other locally finite factor of iid generating set for the cycle space would do the job for this application. In particular, one can replace part c of Theorem 1.2 in \cite{ARS} with the following.

\begin{lemma}\label{generalized}
Let $G$ be a unimodular random graph of finite expected degree at the root and let $x\in [0,1]$. If $\omega$ is a sample from $\FKI$ with $p=\frac{2x}{1+x}$ as in \eqref{fki}, and there is a factor of iid locally finite set $\calC$ of cycles in $\omega$ that generates the cycle space of $\omega$, then the free Loop $O(1)$ model of parameters $x$ and 0 on $G$ is a factor of iid.
\end{lemma}

(We simplified the generality of the claim to restrict ourselves to cases where the question is still open.) 
Cayley graphs of finitely presented groups trivially have factor of iid locally finite generating sets for their cycle space, so an answer to the above question would follow for this class of graphs if one could prove that having such a generating set is preserved by the $\FKI$ clusters. In Theorem \ref{percol} we show this for values of $p$ that are close enough to 1. 

A locally finite factor of iid generating set in the FK-Ising model could be used to generate the Loop $O(1)$ model as a factor of iid. For simplicity, call an invariant locally finite generating set of the cycle space of $G$ a {\it feasible generating set}. When it can be chosen as a factor of iid process, then we call it a factor of iid feasible generating set. 
Does every unimodular random graph have a feasible generating set? In \cite{ARS} the authors suspect that the Gromov monster is a counterexample, but without a proof. 
Is there a simple example of a unimodular random graph with no invariant feasible generating set?
By adding an extra edge to every vertex $v$ of a regular tree, connecting it to a uniformly chosen point of $T$ at distance $r_v$ from $v$, where $r_v$ is some random positive number of infinite expectation and iid over the $v$, then
one can show that in the resulting unimodular random graph for every invariant generating set $\calC$, the expected number of cycles in $\calC$ that contain a random edge of $o$ is infinite.

One may wonder if the property of having a (factor of iid) feasible generating set is preserved under Bernoulli percolation. A preliminary step is Theorem \ref{percol}. We were not trying to find the optimal $\bar p$ given by the proof method. 

\begin{proofof}{Theorem \ref{percol}}
Denote by $G(p)$ the random graph given by Bernoulli($p$) bond percolation on $G$. In the following proof summation is always understood modulo 2.

First we prove that there is an invariant locally finite generating set of cycles for $G(p)$. 
Let $k$ be the length of the longest relator in the given finite presentation of $G$, and set $\calC$ of all cycles of length at most $k$ is a (factor of iid) feasible generating set for $G$. 
We will construct a feasible generating set $\calC'$ for $G(p)$. Every cycle $C$ of $G$ can be written as a mod 2 sum of cycles $C_1,\ldots,C_{k(C)}\in\calC$, and we fix $k(C)$ to be minimal such that there is such a writing. For a set $\calO$ of cycles, $\span (\calO)$ will denote its linear span over the two-element field.

Let $\calC_1'$ be the set of cycles $C\subset G(p)$ with $k(C)=1$ (so $C\in\calC$). Recursively, define
$$\calC_n':=\{C\subset G(p): k(C)=n, C\not\in\span (\calC_{n-1}')\}\cup \calC_{n-1}'.
$$

Then $\calC':=\cup_n \calC_{n}'$ is an invariant generating set for the cycle space of $G(p)$. Let us show that it is locally finite. Suppose by contradiction that some edge $e$ is contained each of the infinitely many distinct cycles $O_1,O_2,\ldots\in\calC'$. Now fix $i$, and cycles $C_j^i\in \calC$ such that $O_i=C_1^i+\ldots+C_{k(O_i)}^i$. (We can observe that $O_i\in\calC_{k(O_i)}'\setminus \calC_{k(O_i)-1}'$.) The subgraph $\cup_j C_j^i$ is connected, otherwise $\sum_j C_j^i$ would have at least two components (also using that $k(O)$ is minimal). Also, at least one of the $C_j^i$ contains $e$. To summarize: $e$ is contained in each of the connected subgraphs $\cup_j C_j^i$ of $G$, and the size of these tends to infinity as $i$ goes to infinity.

Finally, for every $j$ we have $C_j^i\not\subset G(p)$. Because, if we suppose to the contrary that   
$C_j^i\subset G(p)$, then $O_i+C_j^i=\sum_{j'\not=j}C_j^i\subset G(p)$, and this would mean $O_i=C_j+\sum_{j'\not=j}C_j^i\in\span (\calC_{k(O_i)-1}'$, a contradiction. 
We have just seen that every $C_j^i$ has some edge in $G(p)$. Hence every $O_i=\sum_j C_j^i$ is in the $k$-neighborhood ({\it $k$-closure}) of $E(G)\setminus E(G(p))$ in $G$. This $k$-closure therefore has to contain some infinite component. On the other hand the $k$-closure, as a random subgraph of $G$, is stochastically dominated by Bernoulli($q(p)$) bond percolation on $G$ for some small enough $q(p)>0$, and one can express $q(p)$ as a monotone decreasing function of $p$ (decreasing because it needs to dominate a neighborhood of the {\it closed} edges by the $G(p)$ percolation in our setup). See e.g. Proposition 7.14 in \cite{LP}. Take $p<1$ large enough to make $q(p)$ smaller than the critical percolation probability for Bernoulli bond percolation on $G$. This contradicts the existence of infinite components in the $k$-closure of $E(G)\setminus E(G(p))$, finishing the proof that $\calC'$ is a feasible generating set.

To prove that $\calC'$ can be constructed as a factor of iid, we start by observing that whether $k(C)=n$ for a given cycle $C$ can be decided from a bounded neighborhood of $C$, using the fact that $C=\sum_{i=1}^{k(C)} C_i$ implies that $\cup_{i=1}^n C_i$ is connected, as we have seen. Doing it for $n=1,2,\ldots$, we will find $k(C)$ in a finite number of steps, because every $C$ is in $\span (\calC)$. Every cycle $C$ in $\calC_n'$ has length at most $kn$, because it is a sum of at most $n$ elements of $\calC$. So to check if $C\in \calC_n'\setminus \calC_{n-1}$, we first check if $k(C)=n$, $C\subset G(p)$, and then look through all of the finitely many possibilities of writing $C$ as a sum of $n$ elements $C_1,\ldots,C_n$ of $\calC$, using again that $\cup_{i=1}^n C_i$ has to be connected. Then we reduced the problem to checking whether the $C_i$ are in $\calC_{n-1}'$. Examining a bounded neighborhood of $C$ gives us the answer to whether $C\in\calC'$. 

The proof for $\FKI$ works the same way, as before.
\qed
\end{proofof}

\begin{proofof}{Theorem \ref{finitely_pres}}
It follows from Lemma \ref{generalized} and Theorem \ref{percol}.
\qed
\end{proofof}

\bigskip

{\bf Acknowledgments:} 
I thank G\'abor Pete for some comments on the manuscript.

\ \\

{\bf \'Ad\'am Tim\'ar}\\
Division of Mathematics, The Science Institute, University of Iceland\\
Dunhaga 3 IS-107 Reykjavik, Iceland\\
and\\
Alfr\'ed R\'enyi Institute of Mathematics\\
Re\'altanoda u. 13-15, Budapest 1053 Hungary\\
\texttt{madaramit[at]gmail.com}\\


\begin{thebibliography}{AAA}

\bibitem{ADS}
M. Aizenman, H. Duminil-Copin, and V. Sidoravicius (2015) Random currents and
continuity of Ising model’s spontaneous magnetization, {\it Communications in Mathematical
Physics} no. 2, 719–742.

\bibitem{ADTW}
M. Aizenman, H. Duminil-Copin, V. Tassion, and S. Warzel (2019) Emergent
planarity in two-dimensional Ising models with finite-range interactions, {\it Inventiones mathematicae}
216, no. 3, 661–743.

\bibitem{AL}
D. Aldous, R. Lyons (2007) 
Processes on unimodular random networks
{\it Electron. J. Probab.} {\bf 12}, 1454-1508.

\bibitem{ARS}
O. Angel, G. Ray and Y. Spinka (2021) Uniform even subgraphs and graphical
representations of Ising as factors of i.i.d. (preprint)

\bibitem{GJ} G. Grimmett and S. Janson (2007) Random Even Graphs, {\it Electronic Communications in Probability} {\bf 16}.

\bibitem{Cu} N. Curien (2018) Random graphs - the local convergence point of view. Notes \url{https://www.imo.universite-paris-saclay.fr/~nicolas.curien/cours/cours-RG.pdf}

\bibitem{HJL}
O. Häggström, J. Jonasson, and R. Lyons (2002) Coupling and Bernoullicity in random cluster
and Potts models, {\it Bernoulli}, 275–294.





\bibitem{HSr}
M. Harel and Y. Spinka (2022) Finitary codings for the random-cluster model and other
infinite-range monotone models, {\it Electron. J. Probab.} 27: 1-32.
.

\bibitem{H}
T. Hutchcroft (2020) Continuity of the Ising phase transition on nonamenable groups, \url{arXiv:2007.15625}



\bibitem{LP}
R. Lyons and Y. Peres. {\it Probability on Trees and Networks.}
Cambridge University Press, New York, 2016. Available at \url{http://pages.iu.edu/\textasciitilde rdlyons/}



\bibitem{PGG}
G. Pete.
{\it Probability and Geometry on Groups}.
Book in preparation, \url{http://www.math.bme.hu/~gabor/PGG.pdf}



\bibitem{T0} \'A. Tim\'ar (2007) Cutsets in infinite graphs, {\it Combinatorics, Probability and Computing} {\bf 16}, 159-166.


\bibitem{T_fafaktor}
\'A. Tim\'ar (2021) A nonamenable “factor” of a Euclidean space, The Annals of Probability 49, no. 3, 1427–1449.


\end{thebibliography}
\end{document}